\documentclass{amsart}
\usepackage{amssymb}
\usepackage{amsmath}
\usepackage{amsthm}
\usepackage[all]{xy}
\xyoption{2cell}

\newcommand{\V}{{\mathcal V}}

\newcommand{\pu}{{\mathbb P^1}}
\newcommand{\pd}{{\mathbb P^2}}
\newcommand{\proj}{\mathbb P}

\DeclareMathOperator{\loc}{\mathrm{Locus}}
\DeclareMathOperator{\cloc}{\mathrm{ChLocus}}

\newcommand{\ratcurves}{\textrm{Ratcurves}^n(X)}

\newcommand{\cone}{\textrm{NE}}
\DeclareMathOperator{\cycl}{N_1}
\DeclareMathOperator{\pic}{Pic}

\newcommand{\W}{{\mathcal{W}}}
\newcommand{\conx}[1]{\cone\,(#1,X)}
\newcommand{\cycx}[1]{\cycl(#1,X)}

\newcommand{\rc}[2]{#1 \xymatrix{\ar@{-->}[r] & }{#2}}

\newtheorem{theorem}{Theorem}[section]
\newtheorem*{theorem*}{Theorem}
\newtheorem{lemma}[theorem]{Lemma}
\newtheorem{conjecture}[theorem]{Conjecture}
\newtheorem{proposition}[theorem]{Proposition}
\newtheorem{corollary}[theorem]{Corollary}

\theoremstyle{definition}
\newtheorem{definition}[theorem]{Definition}

\theoremstyle{remark}
\newtheorem{remark}[theorem]{Remark}
\newtheorem{construction}[theorem]{Construction}

\begin{document}


\title[Fano manifolds with an unsplit dominating family of curves]{On Fano manifolds with an unsplit dominating family of rational curves}

\renewcommand{\theequation}{{\arabic{section}.\arabic{theorem}.\arabic{equation}}}

\author{Carla Novelli}
\address{Dipartimento di Matematica Pura ed Applicata  \newline
\indent 
Universit\`a degli Studi di Padova \newline
\indent via Trieste, 63  \newline
\indent 
I-35121 Padova, Italy} 
\email{novelli@math.unipd.it}

\begin{abstract}
We study Fano manifolds $X$ admitting an unsplit dominating family of rational curves and we prove that the Generalized Mukai Conjecture holds if $X$ has pseudoindex $i_X = (\dim X)/3$ or dimension $\dim X=6$.
We also show that this conjecture is true for all Fano manifolds with $i_X > (\dim X)/3$. 
\end{abstract}

\maketitle

\section{Introduction}

Let $X$ be a Fano manifold, {\em i.e.} a smooth complex projective variety whose anticanonical bundle $-K_X$ is ample. 
A Fano manifold is associated with two invariants, namely the {\em index}, $r_X$, defined as
the largest integer dividing $-K_X$ in the Picard group of $X$,
and the {\em pseudoindex}, $i_X$, defined as
the minimum anticanonical degree of rational curves on $X$. \par
\medskip

\noindent
In 1988 Mukai proposed the following conjecture, involving the index and the Picard number of a Fano manifold:
\smallskip

\begin{center}
\begin{minipage}[center]{10,5cm}
\begin{conjecture} \cite{Kata} \label{M} \
Let $X$ be a Fano manifold of dimension $n$. \\ Then 
$\rho_X (r_X - 1) \le n$, with equality if and only if $X = (\proj^{r_X-1})^{\rho_X}$.
\end{conjecture}
\end{minipage}
\end{center}\par
\medskip

\noindent
In 1990, in \cite{Wimu}, where the notion of pseudoindex was introduced, the first step towards the conjecture was made and it was proved that if $i_X > ({\dim X+2})/{2}$ then $\rho_X =1$;
moreover, if $r_X = (\dim X+2)/{2}$ then either $\rho_X =1$ or $X = (\proj^{r_X-1})^2$.\par
\medskip

\noindent
In 2002 Bonavero, Casagrande, Debarre and Druel reconsidered this problem and proposed
the following more general conjecture:

\begin{center}
\begin{minipage}[center]{10,5cm}
\begin{conjecture} \cite{BCDD} \label{GM} \
Let $X$ be a Fano manifold of dimension $n$. \\ Then $\rho_X(i_X - 1) \le n$,
with equality if and only if $X = (\proj^{i_X-1})^{\rho_X}$.
\end{conjecture}
\end{minipage}
\end{center}
\noindent
Moreover, in \cite{BCDD}, they proved Conjecture (\ref{GM}) for Fano manifolds of dimension $4$ (in lower dimension the result can be read off from the classification), for homogeneous manifolds, and for toric Fano manifolds of pseudoindex $i_X \geq ({\dim X+3})/{3}$ or dimension $\le 7$.
In 2006, in \cite{Cas06}, the toric case was completely settled.\par
\medskip

\noindent
In 2004, in \cite{ACO}, Conjecture (\ref{GM}) was proved for Fano manifolds of dimension $5$ and for Fano manifolds of pseudoindex $i_X \ge ({\dim X+3})/{3}$ admitting an unsplit dominating family of rational curves (see Definition (\ref{Rf})).\par
\medskip

\noindent
In 2010, in \cite{NO-Mukai}, Conjecture (\ref{GM}) was proved for Fano manifolds of pseudoindex $i_X \ge ({\dim X+3})/{3}$,
and simplified proofs of this conjecture for Fano manifolds of dimension $4$ and $5$ were provided.
\par
\medskip

In this paper we reconsider Fano manifolds $X$ admitting an unsplit dominating family of rational curves, and we prove Conjecture (\ref{GM}) if $X$ has dimension $6$ (Theorem (\ref{six-unsplit})), or $X$ has pseudoindex $i_X \ge (\dim X)/{3}$ (Theorem (\ref{CombiningFDU})).
\par
\medskip

\noindent
The paper is organized as follows: in Sections (\ref{back}) and (\ref{chains}) we recall definitions and results on families of rational curves and on chains of rational curves on projective manifolds, while in Section (\ref{ratcurvFano}) we consider families of rational curves on Fano manifolds;
in Section (\ref{bounding}) we prove Conjecture (\ref{GM}) for Fano manifolds $X$ of pseudoindex $i_X > (\dim X)/{3}$;
in Section (\ref{FanoConFDU}) we consider Fano manifolds $X$ admitting an unsplit dominating family of rational curves and we prove Conjecture (\ref{GM}) if $\dim X=6$, or $i_X \ge (\dim X)/{3}$.

\section{Families of rational curves}\label{back}

Let $X$ be a smooth complex projective variety.\par

\begin{definition} \label{Rf}
A {\em family of rational curves} $V$ on $X$ is an irreducible component
of the scheme $\ratcurves$ (see \cite[Definition II.2.11]{Kob}).\\
Given a rational curve we will call a {\em family of
deformations} of that curve any irreducible component of  $\ratcurves$
containing the point parameterizing that curve.\\
We define $\loc(V)$ to be the set of points of $X$ through which there is a curve among those
parametrized by $V$; we say that $V$ is a {\em covering family} if ${\loc(V)}=X$ and that $V$ is a
{\em dominating family} if $\overline{\loc(V)}=X$.\\
By abuse of notation, given a line bundle $L \in \pic(X)$, we will denote by $L \cdot V$
the intersection number $L \cdot C$, with $C$ any curve among those
parametrized by $V$.\\
We will say that $V$ is {\em unsplit} if it is proper; clearly, an unsplit dominating family is
covering.\\
We denote by $V_x$ the subscheme of $V$ parameterizing rational curves
passing through a point $x$ and by $\loc(V_x)$ the set of points of $X$
through which there is a curve among those parametrized by $V_x$. If, for a general point $x \in \loc(V)$,
$V_x$ is proper, then we will say that the family is {\em locally unsplit}; by Mori's Bend and Break arguments,
if $V$ is a locally unsplit family, then $-K_X \cdot V \le \dim X+1$.\\
If $X$ admits  dominating families, we can choose among them one with minimal degree with respect
to a fixed ample line bundle, and we call it a {\em minimal dominating family}; such
a family is locally unsplit.
\end{definition}

\begin{definition}
Let $U$ be an open dense subset of $X$ and $\pi\colon U \to Z$ a proper
surjective morphism to a quasi-projective variety;
we say that a family of rational curves $V$ is a {\em horizontal dominating family with respect to} $\pi$
if $\loc(V)$ dominates $Z$ and curves parametrized by $V$ are not contracted by $\pi$.
If such families exist, we can choose among them  one with minimal degree with respect
to a fixed ample line bundle and we call it a {\em minimal
horizontal dominating family} with respect to $\pi$; such
a family is locally unsplit.
\end{definition}

\begin{remark} By fundamental results in \cite{Mo79}, a Fano manifold admits dominating families of rational curves;
also horizontal dominating families with respect to proper morphisms defined on an open set exist,
as proved in \cite{KoMiMo}. In the case of Fano manifolds with ``minimal'' we will mean minimal with respect to $-K_X$,
unless otherwise stated.
\end{remark}

\begin{definition}\label{CF}
We define a {\em Chow family of rational 1-cycles} $\W$ to be an irreducible
component of  $\textrm{Chow}(X)$ parameterizing rational and connected 1-cycles.\\
We define $\loc(\W)$ to be the set of points of $X$ through which there is a cycle among those
parametrized by $\W$; notice that $\loc(\W)$ is a closed subset of $X$ (\cite[II.2.3]{Kob}).
We say that $\W$ is a {\em covering family} if $\loc(\W)=X$.\\
If $V$ is a family of rational curves, the closure of the image of
$V$ in $\textrm{Chow}(X)$, denoted by $\V$, is called the {\em Chow family associated to} $V$.
\end{definition}

\begin{remark}
If $V$ is proper, {\em i.e.} if the family is unsplit, then $V$ corresponds to the normalization
of the associated Chow family $\V$.
\end{remark}

\begin{definition}
Let $V$ be a family of rational curves and let $\V$ be the associated Chow family. We say that
$V$ (and also $\V$) is {\em quasi-unsplit} if every component of any reducible cycle parametrized by 
$\V$ has numerical class proportional to the numerical class of a curve parametrized by $V$.
\end{definition}

\begin{definition}
Let $V^1, \dots, V^k$ be families of rational curves on $X$ and $Y \subset X$.\\
We define $\loc(V^1)_Y$ to be the set of points $x \in X$ such that there exists
a curve $C$ among those parametrized by $V^1$ with
$C \cap Y \not = \emptyset$ and $x \in C$. We inductively define
$\loc(V^1, \dots, V^k)_Y := \loc(V^2, \dots, V^k)_{\loc(V^1)_Y}$.
Notice that, by this definition, we have $\loc(V)_x=\loc(V_x)$.
Analogously we define $\loc(\W^1, \dots, \W^k)_Y$  for Chow families $\W^1, \dots, \W^k$ of rational 1-cycles.
\end{definition}

{\bf Notation}: If $\Gamma$ is a $1$-cycle, then we will denote by $[\Gamma]$ its numerical equivalence class
in $\cycl(X)$; if $V$ is a family of rational curves, we will denote by $[V]$ the numerical equivalence class
of any curve among those parametrized by $V$.\\
If $Y \subset X$, we will denote by $\cycx{Y} \subseteq \cycl(X)$ the vector subspace
generated by numerical classes of curves of $X$ contained in $Y$; moreover, we will denote 
by  $\conx{Y} \subseteq \cone(X)$
the subcone generated by numerical classes of curves of $X$ contained in~$Y$.\par
\medskip
We will make frequent use of the following dimensional estimates:

\begin{proposition} (\cite[IV.2.6]{Kob}\label{iowifam})
Let $V$ be a family of rational curves on $X$ and $x \in \loc(V)$ a point such that every component of $V_x$ is proper.
Then
  \begin{itemize}
       \item[(a)] $\dim V-1=\dim \loc(V)+\dim \loc(V_x) \ge \dim X  -K_X \cdot V -1$;
       \item[(b)] $ \dim \loc(V_x) \ge -K_X \cdot V -1$.
    \end{itemize}
\end{proposition}

\begin{definition} We say that $k$ quasi-unsplit families $V^1, \dots, V^k$ are numerically independent
if in $\cycl(X)$ we have $\dim \langle [V^1], \dots, [V^k]\rangle =k$.
\end{definition}

\begin{lemma} (Cf. \cite[Lemma 5.4]{ACO}) \label{locy}
Let $Y \subset X$ be a closed subset and  $V^1, \dots, V^k$ numerically independent
unsplit families of rational curves such that $\langle [V^1], \dots, [V^k]\rangle$ $\cap \conx{Y}={\underline 0}$.
Then either $\loc(V^1, \ldots,V^k)_Y=\emptyset$ or
      $$\dim \loc(V^1, \ldots, V^k)_Y \ge \dim Y +\sum -K_X  \cdot V^i -k.$$
\end{lemma}

A key fact underlying our strategy to obtain bounds on the Picard number, based on \cite[Proposition II.4.19]{Kob},
is the following:

\begin{lemma}\label{numeq} (\cite[Lemma 4.1]{ACO})
Let $Y \subset X$ be a closed subset, $\V$ a Chow family of rational $1$-cycles. Then every curve
contained in $\loc(\V)_Y$ is numerically equivalent to a linear combination with rational
coefficients of a curve contained in $Y$ and of irreducible components of cycles parametrized
by $\V$ which meet $Y$.
\end{lemma}

\begin{corollary}\label{numcor} Let $V^1$ be a locally unsplit family of rational curves,
and $V^2, \dots,$ $V^k$ unsplit families of rational curves. Then, for a general $x \in \loc(V^1)$,
\begin{enumerate}
\item[(a)] $\cycx{\loc(V^1)_x} = \langle [V^1] \rangle$;
\item[(b)] either $\loc(V^1, \dots, V^k)_x= \emptyset$, or 
$\cycx{\loc(V^1, \dots, V^k)_x} = \langle [V^1], \dots,$ $[V^k] \rangle$.
\end{enumerate}
\end{corollary}

\section{Chains of rational curves}\label{chains}

Let $X$ be a smooth complex projective variety.
Let $V$ be a dominating family of rational curves on $X$ and denote by $\V$
the associated Chow family.

\begin{definition}
Let $Y \subset X$ be a closed subset; define $\cloc_m(\V)_Y$
to be the set of points $x \in X$ such that there exist cycles $\Gamma_1, \dots, \Gamma_m$
with the following properties:
    \begin{itemize}
       \item $\Gamma_i$ belongs to the family $\V$;
       \item $\Gamma_i \cap \Gamma_{i+1} \not = \emptyset$;
       \item $\Gamma_1 \cap Y \not = \emptyset$ and $x \in \Gamma_m$,
    \end{itemize}
{\em i.e.} $\cloc_m(\V)_Y$ is the set of points that can be
joined to $Y$ by a connected chain  of at most $m$ cycles belonging to the family $\V$.\\
If we consider among cycles parametrized by $\V$ only irreducible ones, 
in the same way we can define $\cloc_m(V)_Y$.

\end{definition}

Define a relation of {\em rational connectedness with respect to
$\V$} on $X$ in the following way: two points $x$ and $y$ of $X$ are in rc$(\V)$-relation if
there exists a chain of cycles in $\V$ which joins
$x$ and $y$, {\em i.e.} if $y \in \cloc_m(\V)_x$ for some $m$.
In particular,  $X$ is {\em $rc(\V)$-connected} if for some $m$ we have
$X=\cloc_m(\V)_x$.\par
\medskip
The family $\V$ defines a proper prerelation in the sense of \cite[Definition IV.4.6]{Kob}.
This prerelation is associated with a fibration, which we will call the {\em rc$(\V)$-fibration}:

\begin{theorem}(\cite[IV.4.16]{Kob}, Cf. \cite{Cam81}) Let $X$ be a normal and proper variety and $\V$
a proper prerelation; then there exists an open subvariety $X^0 \subset X$ and a proper morphism with
connected fibers $\pi\colon X^0 \to Z$ such that
\begin{itemize}
\item  $\langle \mspace{1mu} {\mathcal U}\mspace{1mu} \rangle$ restricts to an equivalence relation on $X^0$;
\item $\pi^{-1}(z)$ is a  $\langle \mspace{1mu} {\mathcal U} \mspace{1mu}\rangle$-equivalence class for every $z \in Z$;
\item $\forall\, z\in Z$ and  $\forall\, x,y \in \pi^{-1}(z)$, $x \in \cloc_m(\V)_y$ with
$m \le 2^{\dim X -\dim Z}-1$.
\end{itemize}
\end{theorem}

Clearly $X$ is rc$(\V)$-connected if and only if $\dim Z^0=0$.

\medskip

Given $\V^1, \dots, \V^k$ Chow families of rational 1-cycles, it is possible to define a relation of
rc$(\V^1, \dots, \V^k)$-connectedness, which is associated with a fibration, that we will call
rc$(\V^1, \dots, \V^k)$-fibration. The variety $X$ will be called  {\em rc$(\V^1, \dots, \V^k)$-connected} if
the target of the fibration is a point.\par
\medskip
For such varieties we have the following application of Lemma (\ref{numeq}):
\begin{proposition} (Cf. \cite[Corollary 4.4]{ACO})\label{rhobound}
If $X$ is rationally connected with respect to some Chow families of rational $1$-cycles
$\V^1, \dots, \V^k$, then $\cycl(X)$ is generated by the classes of irreducible components of cycles
in $\V^1, \dots, \V^k$.\\
In particular, if $\V^1, \dots, \V^k$ are quasi-unsplit families, then $\rho_X \le k$ and equality
holds if and only if $\V^1, \dots, \V^k$ are numerically independent.
\end{proposition}

A straightforward consequence of the above proposition is the following:

\begin{corollary}\label{rhoboundcor}(\cite[Corollary 3]{NO-Mukai})
If $X$ is rationally connected with respect to Chow families of rational $1$-cycles
$\V^1, \dots, \V^k$ and $D$ is an effective divisor, then $D$ cannot be trivial on every irreducible component
of every cycle parametrized by $\V^1, \dots, \V^k$.
\end{corollary}

We will also make use of the following

\begin{lemma}\label{freq}(\cite[Lemma 3]{NO-Mukai}) Let $X$ be a Fano manifold of pseudoindex $i_X$, let $Y \subset X$ be a closed subset
of dimension $\dim Y > \dim X - i_X$
and let $W$ be an unsplit non dominating family of rational curves such that $[W] \not \in \conx{Y}$.
Then $\loc(W) \cap Y = \emptyset$.
\end{lemma}

\section{Families of rational curves on Fano manifolds}\label{ratcurvFano}

We start this section by recalling the following

\begin{construction}\label{kfam}(\cite[Construction 1]{NO-Mukai}) Let $X$ be a Fano manifold; let
$V^1$ be a minimal dominating family of rational curves on $X$ and consider the associated Chow family $\V^1$.\\
If $X$ is not rc$(\V^1)$-connected, let $V^2$ be a minimal
horizontal dominating family with respect to the rc$(\V^1)$-fibration,
$\pi_{1}\colon X^{} \xymatrix{\ar@{-->}[r]^{}&} Z^{1}$. 
If $X$ is not rc$(\V^1, \V^2)$-connected,
we denote by $V^3$ a minimal horizontal dominating family with respect to the
the rc$(\V^1,\V^2)$-fibration, $\pi_{2}\colon X^{} \xymatrix{\ar@{-->}[r]^{}&} Z^{2}$, and so on.
Since $\dim Z^{i+1} < \dim Z^i$, for some integer $k$ we have that $X$ is
rc$(\V^1, \dots, \V^k)$-connected.\\
Notice that, by construction, the families $V^1, \dots, V^k$ are numerically independent.
\end{construction}

\begin{lemma}\label{kfamprop}(\cite[Lemma 4]{NO-Mukai})
Let $X$ be a Fano manifold of pseudoindex $i_X \ge 2$ and let $V^1, \dots, V^k$ be families
of rational curves as in Construction (\ref{kfam}). Then
$$\sum_{i=1}^k(-K_X \cdot V^i -1) \le \dim X.$$
In particular, $k (i_X -1) \le \dim X$, and equality holds if and only if $X = (\proj^{i_X-1})^k$.
\end{lemma}

\begin{lemma}\label{k-2fam}
Let $X$ be a Fano manifold of pseudoindex $i_X \ge 2$ and let $V^1, \dots, V^k$ be families
of rational curves as in Construction (\ref{kfam}). 
Assume that at least one of these families, say $V^j$, is not unsplit. Then \; $k (i_X-1)\leq \dim X - i_X$.

Moreover, 
\begin{itemize}
	\item[(a)] if $j = \frac{\dim X - i_X}{i_X-1}$, then $j=k$ and \; $\rho_X (i_X-1)=\dim X - i_X$;
	\item[(b)] if $j = \frac{\dim X - i_X-1}{i_X-1}$, then $j=k$ and either \; $\rho_X (i_X-1)=\dim X - i_X -1$, or \; $i_X=2$ and $\rho_X=\dim X -2$.
\end{itemize}
\end{lemma}

\proof
Let $V^1, \dots, V^k$ be families of rational curves as in Construction (\ref{kfam});
by Lemma (\ref{kfamprop}) we get $(k-1)(i_X-1)+(2i_X-1)\le \dim X$, hence $k \le \frac{\dim X - i_X}{i_X-1}$.
Moreover, by part (b) of Proposition (\ref{iowifam}), we have $\dim \loc(V^j)_{x_j} \ge 2i_X-1$ for a general point $x_j \in \loc(V^j)$.

\medskip

If $j=\frac{\dim X - i_X}{i_X-1}$, then $j=k$ and $V^j$ is the only non unsplit family.
Then, for a general point $x_k \in \loc (V^k)$, we have 
$X= \loc(V^k, \dots ,V^1)_{x_k}$ by Lemma (\ref{locy}).
Therefore, by part (b) of Corollary (\ref{numcor}), we obtain that
$\cycl(X) = \langle [V^1], \dots, [V^k]\rangle$, so $\rho_X=k$, and we obtain case (a) of the statement.

\smallskip

Assume now that $j=\frac{\dim X - i_X-1}{i_X-1}$. Then $V^j$ is the only non unsplit family; moreover, $\dim$ $\loc(V^j,\dots,V^1)_{x_j} \geq\dim X-1$ by Lemma (\ref{locy}).

We claim that $X$ is rc$(V^1,\dots, \V^j)$-connected.\\
In fact, a general fiber of the rc$(V^1,\dots, \V^j)$-fibration has dimension at least $\dim$ $\loc(V^j,\dots,V^1)_{x_j}\ge \dim X-1$ by Lemma (\ref{locy}).
This implies $\dim Z^j \le 1$, and thus, if $X$ were not rc$(V^1,\dots,$ $\V^j)$-connected, we would have $\dim \loc(V^{j+1})_{x_{j+1}}=1$ for a general point $x_{j+1} \in \loc(V^{j+1})$.
Hence, by part (b) of Proposition (\ref{iowifam}), $-K_X \cdot V^{j+1} =2=i_X$, so $V^{j+1}$ would be unsplit and, by part (a) of the same proposition, covering, against the minimality of $V^j$.
Therefore $j=k$.

Consider an irreducible component $D$ of $\loc(V^k,\dots,V^1)_{x_k}$ of maximal dimension (which is at least $\dim X -1$). Therefore, either $X=\loc(V^k,\dots,V^1)_{x_k}$ and $\rho_X=k$ by part (b) of Corollary (\ref{numcor}),
or $D$ is a divisor in $X$. In this last case, $N_1(D,X)=\langle [V^1], \dots, [V^k] \rangle$ by part (b) of Corollary (\ref{numcor}).
Then, by \cite[Theorem 1.6 and Corollary 2.12]{Cas09}, either $\rho_X=k$, or $i_X=2$ and $\rho_X=\dim X-2$.
\qed

\section{Bounds on the Picard number of Fano manifolds}\label{bounding}

In this section we show that Conjecture (\ref{GM}) holds for Fano manifolds of pseudoindex $i_X > \dim X/{3}$.

\begin{theorem}\label{terzini}
Let $X$ be a Fano manifold of Picard number $\rho_X$ and pseudoindex $i_X > \dim X/{3}$. Then
$\rho_X (i_X-1) \leq \dim X$ and equality holds if and only if $X=(\proj^{i_X-1})^{\rho_X}$.
\end{theorem}

\begin{proof}
Note that in view of \cite[Theorem 3]{NO-Mukai} we can restrict to $i_X < (\dim X+3)/{3}$.
Moreover, since for $i_X=1$ there is nothing to prove, we assume $i_X\ge 2$ (and so $\dim X>3$).\par
\smallskip

Let $V^1, \dots, V^k$ be families of rational curves as in Construction (\ref{kfam}).\par

If all the families are unsplit, then Lemma (\ref{kfamprop}) gives $k \le 3$ unless either 
$i_X=2$, $\dim X=5$ and $k=4$, or $X=(\pu)^5$, or $X=(\pu)^4$, or $X=(\pd)^4$.
\\
Since $\rho_X =k$ by Proposition (\ref{rhobound}), the assertion follows.\par
\smallskip

We can thus assume that at least one of these families, say $V^j$, is not unsplit. 
Then, by Lemma (\ref{kfamprop}), $k \le 3$ and exactly one of these families is not unsplit.
Moreover, if $j=3$, by computing $\dim\loc(V^3,V^2,V^1)$ with Lemma (\ref{locy}), we get a contradiction unless $\dim X=5$ and $i_X=2$, so $\rho_X=3$ by part (b) of Corollary (\ref{numcor}).
If $j=2$ and $i_X=(\dim X+2)/3$, then $\rho_X=2$ by part (a) of Lemma (\ref{k-2fam}).
If $j=2$ and $i_X=(\dim X+1)/3$, denoted by $T$ an irreducible component of maximal dimension of $\loc(V^2,V^1)_{x_2}$, we have $\dim T\geq \dim X-1$ by Lemma (\ref{locy}).
Since $\cycl(T,X)= \langle [V^1], [V^2] \rangle$ by part (b) of Corollary (\ref{numcor}), we have that if $\dim T=\dim X$ then $\rho_X=2$, while if $\dim T=\dim X-1$ then either $\rho_X=2$ or $\dim X=5$, $i_X=2$ and $\rho_X=3$ by \cite[Theorem 1.6 and Corollary 2.12]{Cas09}.
\par
\smallskip

Therefore we are left with $j=1$.
Then a general fiber of the rc$(\V^1)$-fibration $X^{} \xymatrix{\ar@{-->}[r]^{}&} Z^{1}$ has dimension at least $\dim \loc(V^1)_{x_1}$.

Assume first that $\dim Z^1 \ge 1$. 
Since for a general point $x_2 \in \loc(V^2)$ we know that $\dim \loc(V^2)_{x_2} \leq \dim Z^1$, we deduce that $-K_X\cdot V^2\leq i_X+1$ by part (b) of Proposition (\ref{iowifam}). So $V^2$ is unsplit and $V^2$ is not dominating, since $-K_X\cdot V^2 < -K_X\cdot V^1$.
Denote by $D$ an irreducible component of maximal dimension of $\loc(V^1,V^2)_{x_1}$. Then $\dim D=\dim X-1$ and $\cycl(D,X)= \langle [V^1], [V^2] \rangle$, so we are done by \cite[Theorem 1.6 and Corollary 2.12]{Cas09}.

Finally we deal with the case in which $\dim Z^1 =0$, so $X$ is rc$(\V^1)$-connected.
Let $x$ be a general point. 
Since $x$ is general and $V^1$ is minimal we have $\overline{\loc(V^1)_x}=\loc(V^1)_x$ and 
$\cycx{\loc(V^1)_x}=\langle [V^1] \rangle$ by part (a) of Corollary (\ref{numcor}).\\
If $\loc(V^1)_x=X$, then  $\rho_X =1$. 
So we can suppose that $\dim \loc(V^1)_x < \dim X$ and thus, by part (b) of Proposition (\ref{iowifam}), $-K_X \cdot V^1 \leq \dim X$. 
In particular every reducible cycle parametrized by $\V^1$ has at most two irreducible components.\\
If every irreducible component of a $\V^1$-cycle in a connected $m$-chain through $x$ is numerically proportional to $V^1$,
then $\rho_X=1$ by repeated applications of Lemma (\ref{numeq}).
\\
We can thus assume that there exist $m$-chains through $x$, $\Gamma_1 \cup \Gamma_2 \cup \dots \cup \Gamma_m$, 
with $x \in \Gamma_1$ and $\Gamma_i \cap \Gamma_{i+1} \not = \emptyset$, such that, 
for some $j \in \{1, \dots, m\}$ the irreducible components $\Gamma_j^1$ and $\Gamma_j^2$ of $\Gamma_j$
are not numerically proportional to $\V^1$.\\
Let $j_0 \in \{1, \dots, m\}$ be the minimum integer for which such a chain exists; by the generality of $x$ we have $j_0 \ge 2$.
If $j_0=2$ set $x_1=x$, otherwise let $x_1$ be a point in $\Gamma_{j_0-1} \cap \Gamma_{j_0 -2}$. 
Since $\Gamma_{j_0-1} \subset \loc (\V^1)_{x_1}$ there is  an irreducible component $Y$ of $\loc(V^1)_{x_1}$ which meets $\Gamma_{j_0}$. 
By Lemma (\ref{numeq}), $\cycx{Y} = \langle [V^1] \rangle$.\\
Let $\gamma$ be a component of $\Gamma_{j_0}$ meeting $Y$ and denote by $W$ a family of deformations of $\gamma$; 
then the family $W$ is unsplit and it is not covering, by the minimality of $V^1$.\\
Then $\dim \loc(W)_{Y} =\dim X-1$, and so $\loc(W)=\loc(W)_{Y}$. Moreover, in this case,
by part (b) of Corollary (\ref{numcor}) we get $\cycl(\loc(W)_{Y},X)= \langle [V^1], [W] \rangle$.\\
Therefore $\rho_X=2$ by \cite[Theorem 1.6 and Corollary 2.12]{Cas09}.
\end{proof}

\smallskip

Now, in view of Theorem (\ref{terzini}) it is straightforward to derive the following results.

\begin{proposition}\label{fraz}
Let $X$ be a Fano manifold of dimension $\geq 7$, Picard number $\rho_X$ and pseudoindex $i_X > (\dim X-3)/{2}$. Then
$\rho_X (i_X-1) \leq \dim X$ and equality holds if and only if $X=(\proj^{i_X-1})^{\rho_X}$.
\end{proposition}

\begin{proposition}\label{-4}
Let $X$ be a Fano manifold of Picard number $\rho_X$ and pseudoindex $i_X > \dim X -4$. Then
$\rho_X (i_X-1) \leq \dim X$ and equality holds if and only if $X=(\proj^{i_X-1})^{\rho_X}$.
\end{proposition}

\begin{remark}
All the previous results can be improved once the Generalized Mukai Conjecture is proved in the case of Fano manifolds of dimension $6$. However, this seems to be much more difficult, so in the next section we prove the conjecture under some additional assumption.
\end{remark}

\section{Fano manifolds with an unsplit dominating family}\label{FanoConFDU}
Since the Generalized Mukai Conjecture holds for Fano manifolds of dimension lower than or equal to five, in the next theorems we deal with manifolds of dimension at least six: in Theorem (\ref {>six-unsplit}) we consider Fano manifolds of dimension grater than six and pseudoindex $\dim X / 3$, while in Theorem (\ref {six-unsplit}) we consider Fano sixfolds.\par
\smallskip

We start with the following

\begin{lemma}\label{lemma-terzi}
Let $X$ be a Fano manifold of Picard number $\rho_X$ and pseudoindex $i_X=\dim X /3$. 
If $X$ admits an unsplit dominating family $V$ of rational curves such that $-K_X\cdot V > \dim X /3$, then
$\rho_X (i_X-1) < \dim X $.
\end{lemma}

\begin{proof}
Note that for $i_X=1$ there is nothing to prove, so we can assume $i_X\ge 2$ (and so $\dim X\geq 6$).

Since $V$ is an unsplit dominating family of rational curves on $X$, then either $X$ is rc$(V)$-connected and so $\rho_X=1$, or there exists a minimal horizontal dominating family $V'$ with respect to the rc$(V)$-fibration.\\
In this last case, if $V'$ is not unsplit, we get that an irreducible component $D$ of $\loc(V',V)_{x'}$, for a general point $x' \in \loc(V')$, has dimension at least $\dim X-1$ by Lemma (\ref{locy}).
By part (b) of Corollary (\ref{numcor}), $N_1(D,X)=\langle [V], [V']\rangle$, so, by \cite[Theorem 1.6 and Corollary~2.12]{Cas09}, we have $\rho_X=2$ unless $\dim X=6$ and $\rho_X=3$.\\
We can thus assume that $V'$ is unplit. Now, either $X$ is rc$(V,V')$-connected and so $\rho_X=2$, or there exists a minimal horizontal dominating family $V''$ with respect to the rc$(V,V')$-fibration.
If $V''$ is not unsplit, then by Lemma (\ref{locy}) we can compute $\dim\loc(V'',V',V)_{x''}$ for a general point $x'' \in \loc(V'')$; then we reach a contradiction unless $\dim X=6$ and $\rho_X=3$ by part (b) of Corollary (\ref{numcor}).
If otherwise $V''$ is unsplit, then either $X$ is rc$(V,V',V'')$-connected and so $\rho_X=3$, or there exists a minimal horizontal dominating family $V'''$ with respect to the rc$(V,V',V'')$-fibration.
Then, for a general point $x''' \in \loc(V''')$, computing the dimension of $\loc(V''',V'',V',V)_{x'''}$, we find that either $\dim X=6$ or $9$, $X=\loc(V''',V'',V',V)_{x'''}$ and $\rho_X=4$, or $\dim X=6$, an irreducible component of maximal dimension of $\loc(V''',V'',V',V)_{x'''}$ is a divisor and $\rho_X=4$, or $5$, by part (b) of Corollary (\ref{numcor}) and by \cite[Theorem 1.6 and Corollary 2.12]{Cas09} since $N_1(D,X)=\langle [V], [V'], [V''], [V''']\rangle$.
\end{proof}

\begin{theorem}\label{>six-unsplit}
Let $X$ be a Fano manifold of Picard number $\rho_X$, dimension $\dim X >6$ and pseudoindex $i_X=\dim X /3$. 
If $X$ admits an unsplit dominating family of rational curves, then
$\rho_X (i_X-1) \leq \dim X$ and equality holds if and only if $X=(\proj^3)^4$.
\end{theorem}

\begin{proof}
Denote by $V$ any unsplit dominating family of rational curves on $X$.
We can assume that $-K_X\cdot V=\dim X/3$, since if there exists an unsplit dominating family such that $-K_X \cdot V > \dim X/3$, the assertion follows by Lemma (\ref{lemma-terzi}).
Let $V^1, \dots, V^k$ be families of rational curves as in Construction (\ref{kfam}); then
by Lemma (\ref{kfamprop}) we get $k \le 3$, unless $k=4$, $\dim X=9$ and $i_X=3$, or $X=(\proj^3)^4$.\par

If all the families $V^i$ are unsplit, then $\rho_X =k$ by Proposition (\ref{rhobound}).
\par

We can thus assume that at least one of these families, say $V^j$, is not unsplit. 
Since $-K_X \cdot V^j \ge 2 \dim X /3$, by Lemma (\ref{kfamprop}) we can have only one non-unsplit family among $V^2, \dots, V^k$ and $k \leq 3$.
Moreover, if $j=3$, then $\dim X=9$ by Lemma (\ref{kfamprop}), so $\rho_X =3$ by part (a) of Lemma (\ref{k-2fam}).
\\
So we are left to consider $j=2$. We claim that in this case $X$ is rc$(V^1,\V^2)$-connected. In fact, if this were not the case, there should be a family $V^3$ which is horizontal with respect to the rc$(V^1,\V^2)$-fibration. 
Then, by Lemma (\ref{kfamprop}), we would have that $\dim X=9$ and, by Proposition (\ref{iowifam}), that all the $V^i$s are dominating with $-K_X\cdot V^2 > -K_X\cdot V^3$, which is a contradiction. 
\\
Consider an irreducible component $G$ of $\loc(V^2,V^1)_{x_2}$ of maximal dimension. 
Then $\dim G \ge \dim X -2$ by Lemma (\ref{locy}) and $N_1(G,X)=\langle [V^1],[V^2] \rangle$ by part (b) of Corollary (\ref{numcor}).
If $\dim G = \dim X$ then clearly $\rho_X=2$, while if $\dim G = \dim X-1$ then $\rho_X=2$ by \cite[Theorem 1.6 and Corollary 2.12]{Cas09}.\\
We can thus assume that $\dim G = \dim X-2$.
Since, if all the components of these cycles are contained in $\langle [V^1], [V^2] \rangle$ then $\rho_X=2$,
we can assume that this is not the case. 
Let $\Gamma =\Gamma_1 + \Gamma_2$ be a reducible cycle of $\V^2$ which is not contained in $\langle [V^1], [V^2] \rangle$ and denote by $W^i$ a family of deformations of $\Gamma_i$, $i=1,2$.\\
By Lemma (\ref{locy}) we get $-K_X \cdot V^1=i_X$, $-K_X \cdot V^2=2i_X$ and $\dim \loc(V^2)_{x_2}=2i_X-1$, so that $V^2$ is covering by Proposition (\ref{iowifam}).
\\
We claim that there does not exist any $W^i$, among the families that are not contained in $\langle [V^1], [V^2] \rangle$, such that  $\dim \loc(W^i)=\dim X-1$. In fact, if such a family $W^i$ exist, then it could not be trivial on both $V^1$ and $V^2$ by Corollary (\ref{rhoboundcor}) and Lemma (\ref{freq}).
Therefore $\loc(W^i)$ would intersect $\loc(V^2,V^1)_{x_2}$, so $\dim \loc(V^2,V^1,W^i)_{x_2}$ $\geq \dim X$, which is a contradiction since $W^i$ is not covering.
\\
It follows that $\dim \loc(W^i)\leq \dim X-2$ for any family $W^i$ that is not contained in $\langle [V^1], [V^2] \rangle$.
Then $\loc(W^1,W^2,V^1)_x$ has an irreducible component $D$ of dimension at least $\dim X-1$ by combining Lemma (\ref{locy}) and part (b) of Proposition (\ref{iowifam}).
As $N_1(D,X)=\langle [V^1], [W^1],[W^2] \rangle$ by part (b) of Corollary (\ref{numcor}), we conclude that $\rho_X=3$: this is clear if $\dim D=\dim X$, while it follows by \cite[Theorem 1.6 and Corollary~2.12]{Cas09} if $\dim D=\dim X-1$.
\end{proof}

\begin{theorem}\label{six-unsplit}
Let $X$ be a Fano manifold of Picard number $\rho_X$, pseudoindex $i_X$ and dimension $6$. 
If $X$ admits an unsplit dominating family of rational curves, then
$\rho_X (i_X-1) \leq 6$.
Moreover, equality holds if and only if $X=\proj^6$, or $X=\proj^3 \times \proj^3$, or $X=\proj^2 \times \proj^2\times \proj^2$, or
$X=\proj^1 \times \proj^1 \times \proj^1 \times \proj^1 \times \proj^1 \times \proj^1$.
\end{theorem}

\begin{proof} Clearly we can assume $i_X \geq 2$. Moreover, we can restrict to $i_X = 2$, since otherwise we can apply \cite[Theorem 3]{NO-Mukai}.
So we have to show that $\rho_X \leq 6$, with equality if and only if $X=(\proj^1)^6$.\par

Denote by $V$ any unsplit dominating family of rational curves on $X$.
We can assume that $-K_X\cdot V=2$, since if there exists an unsplit dominating family such that $-K_X\cdot V\geq 3$ then the assertion follows by Lemma (\ref{lemma-terzi}).
Let $V^1, \dots, V^k$ be families of rational curves as in Construction (\ref{kfam}); then
by Lemma (\ref{kfamprop}) we get $k \le 5$, unless $X=(\proj^1)^6$.\par
\smallskip

If all the families $V^i$ are unsplit, then $\rho_X =k$ by Proposition (\ref{rhobound}).
\smallskip

We can thus assume that at least one of these families, say $V^j$, is not unsplit. 
Since $-K_X \cdot V^j \ge 4$, by Lemma (\ref{kfamprop}) we can have only one non-unsplit family among $V^2, \dots, V^k$ and $k \leq 4$.
Moreover, if $j=4$, then $\rho_X =4$ by part (a) of Lemma (\ref{k-2fam}),
while, if $j=3$, then we conclude by part (b) of the same lemma. \par
\smallskip

Therefore we are left with $j=2$. In this case, a general fiber of the rc$(V^1,\V^2)$-fibration $\pi_{2}\colon X^{} \xymatrix{\ar@{-->}[r]^{}&} Z^{2}$ has dimension at least $\dim \loc(V^2,V^1)_{x_2}$, which is at least four by combining Lemma (\ref{locy}) and part (b) of Proposition (\ref{iowifam}).
Then $\dim Z^2 \le 2$.\par

Assume first that $\dim Z^2 \ge 1$ and denote by $V^3$ a minimal horizontal dominating family with respect to $\pi_2$.
Then $\dim \loc(V^3)_{x_3} \leq 2$, so $-K_X \cdot V^3 \leq 3$, by part (b) of Proposition (\ref{iowifam}), and $V^3$ is unsplit.
Moreover, if  $-K_X \cdot V^3 = 3$, then $V^3$ would be covering by Proposition (\ref{iowifam}), contradicting the minimality of $V^2$.
Therefore $-K_X \cdot V^3 = 2$; since $V^3$ cannot be covering, the same proposition implies that $\dim \loc(V^3)_{x_3} = 2$.
It follows that $X$ is rc$(V^1,\V^2,V^3)$-connected.
\\
We claim that $\rho_X=3$.
Let $F$ be a general fiber of the rc$(V^1,\V^2)$-fibration, whose dimension is equal to four. 
Consider an irreducible component $D$ of $\loc(V^3)_F$ of maximal dimension.
By Lemma (\ref{locy}),
$D$ is a divisor.
If $D\cdot V^1 >0$, then, being $V^1$ covering, $X=\loc(V^1)_D$, and $\rho_X=3$ by part (b) of Corollary (\ref{numcor}).
Therefore we can assume $D\cdot V^1 =0$. 
Moreover, if all the components of all the reducible cycles of $\V^2$ are contained in $\langle [V^1], [V^2], [V^3] \rangle$, then $\rho_X=3$,
so we can assume that this is not the case. 
Let $\Gamma =\Gamma_1 + \Gamma_2$ be a reducible cycle of $\V^2$ not contained in $\langle [V^1], [V^2], [V^3] \rangle$.
Then by Lemma (\ref{freq}) $D\cdot \Gamma_i=0$, for $i=1,2$. So $D\cdot V^2 =0$, hence 
we get a contradiction since $V^3$ cannot be trivial on both $V^1$ and $V^2$.

Assume now that $\dim Z^2 =0$, so that $X$ is rc$(V^1, \V^2)$-connected.

If $-K_X \cdot V^2 \geq 6$, then Lemma (\ref{kfamprop}) implies that $-K_X \cdot V^2 = 6$.
It follows by Lemma (\ref{locy}) that $X=\loc(V^2,V^1)_{x_2}$, for a general $x_2 \in \loc(V^2)$ and $\rho_X=2$ by part (b) of Corollary (\ref{numcor}).

Therefore we can assume that $-K_X \cdot V^2 < 6$, so that the reducible cycles of $\V^2$ have exacly two irreducible components.
Consider an irreducible component $G$ of $\loc(V^2,V^1)_{x_2}$ of maximal dimension. Then $\dim G \ge 4$ by Lemma (\ref{locy}).\\
Moreover, if $\dim G=6$, then $\rho_X=2$, so we need to consider $\dim G=4$ or~$5$.
Since, if all the components of these cycles are contained in $\langle [V^1], [V^2] \rangle$ then $\rho_X=2$,
we can assume that this is not the case. 
Let $\Gamma =\Gamma_1 + \Gamma_2$ be a reducible cycle of $\V^2$ not contained in $\langle [V^1], [V^2] \rangle$ and denote by $W^i$ a family of deformations of $\Gamma_i$, $i=1,2$.\\
If $\dim G=5$, then by Lemma (\ref{freq}) $G\cdot \Gamma_i=0$, for $i=1,2$. It follows that $G\cdot V^2 =0$, whence $G\cdot V^1 >0$ by Corollary~(\ref{rhoboundcor}).
Then $X=\loc(V^1)_G$, so $\cycl(X)=\langle [V^1], [V^2] \rangle$, a contradiction.\\
Therefore we are left with $\dim G=4$. By Proposition (\ref{iowifam}) we get $-K_X \cdot V^1=2$, $-K_X \cdot V^2=4$ and $\dim \loc(V^2)_{x_2}=3$, so $V^2$ is covering.
\\
Assume that there exists a family $W^i$, among the families that are not contained in $\langle [V^1], [V^2] \rangle$, such that  $\dim \loc(W^i)=5$. Then it cannot be trivial on both $V^1$ and $V^2$ by Corollary (\ref{rhoboundcor}) and Lemma (\ref{freq}).
Therefore $\loc(W^i)$ intersects $\loc(V^2,V^1)_{x_2}$, so $\dim \loc(V^2,V^1,$ $W^i)_{x_2}$ $=5$, so we conclude by \cite[Theorem 1.6 and Corollary 2.12]{Cas09}.
\\
We can thus assume that $\dim \loc(W^i)=4$ for any family that is not contained in $\langle [V^1], [V^2] \rangle$.
Then $\loc(W^1,W^2,V^1)_{y_1}$ has an irreducible component $D$ of dimension at least five by Lemma (\ref{locy}).
Since $N_1(D,X)=\langle [V^1], [W^1],[W^2] \rangle$, we conclude by part (b) of Corollary (\ref{numcor}) if $\dim D=6$ and by \cite[Theorem 1.6 and Corollary 2.12]{Cas09} if $\dim D=5$.
\end{proof}

By combining the results of this section we actually have the following

\begin{theorem}\label{CombiningFDU}
Let $X$ be a Fano manifold of Picard number $\rho_X$ and pseudoindex $i_X \geq \mbox{\rm min} \{\dim X-4, \ (\dim X-2)/{2}, \ \dim X /3\}$. If $X$ admits an unsplit dominating family of rational curves, then
$\rho_X (i_X-1) \leq \dim X$ and equality holds if and only if $X=(\proj^{i_X-1})^{\rho_X}$.
\end{theorem}

\bigskip

\noindent
\small{{\bf Acknowledgements}. 
I would like to thank Gianluca Occhetta for helpful suggestions during the preparation of the paper.

\end{document}